\documentclass[12pt,reqno]{amsart}

\usepackage{etex}
\usepackage[utf8]{inputenc}
\usepackage{amsfonts}
\usepackage{amsmath}
\usepackage{amssymb}
\usepackage{amsthm}
\usepackage{mathrsfs}
\usepackage{stmaryrd}
\usepackage{color}
\usepackage[english]{babel}
\usepackage[T1]{fontenc}
\usepackage{url}
\usepackage{graphicx}
\usepackage[pdftex,colorlinks,backref=page,citecolor=blue]{hyperref}
\usepackage{caption}
\usepackage{enumerate}
\usepackage{epstopdf}
\usepackage{hyphenat}
\usepackage{float}
\usepackage{indentfirst}
\usepackage[export]{adjustbox}
\usepackage{tikz}
\usepackage[all]{xy}
\usetikzlibrary{matrix,arrows,decorations.pathmorphing}
\usepackage{booktabs}
\usepackage{array}
\usepackage{bm}
\usepackage{csquotes}
\usepackage{enumerate}
\usepackage{comment}
\usepackage{multicol}

\newcolumntype{P}[1]{>{\centering\arraybackslash}p{#1}}

\allowdisplaybreaks

\newtheorem{theorem}{Theorem}[section]

\newtheorem{lemma}[theorem]{Lemma}
\newtheorem{proposition}[theorem]{Proposition}
\newtheorem{corollary}[theorem]{Corollary}

\theoremstyle{definition}

\newtheorem{remark}[theorem]{Remark}
\newtheorem{notation}[theorem]{Notation}
\newtheorem{example}[theorem]{Example}
\newtheorem{definition}[theorem]{Definition}
\newtheorem{question}[theorem]{Question}

\newcommand{\C}{\mathbb{C}}

\newcommand{\N}{\mathbb{N}}
\newcommand{\Q}{\mathbb{Q}}
\newcommand{\KK}{\Bbbk}
\newcommand{\eHK}{e_{\mathrm{HK}}}
\newcommand{\HK}{\mathrm{HK}}
\newcommand{\HKS}{\mathrm{HKS}}
\newcommand{\FS}{\mathrm{FS}}
\newcommand{\FSS}{\mathrm{FSS}}
\newcommand{\fs}{\mathrm{s}}

\DeclareMathOperator{\frk}{frk}

\makeatletter
\def\l@subsection{\@tocline{2}{0pt}{2.5pc}{5pc}{}}
\makeatother

\author{Alessio Caminata}
\address{Alessio Caminata, Dipartimento di Matematica, Dipartimento di Eccellenza 2023-2027, Universit\`a di Genova\\ via Dodecaneso 35, 16146, Genova, Italy}
\email{alessio.caminata@unige.it}

\author{Francesco Zerman}
\address{Francesco Zerman, UniDistance Suisse, Schinerstrasse 18, 3900, Brig, Switzerland}
\email{francesco.zerman@unidistance.ch}

\title[Hilbert-Kunz series, F-signature series, and weak $p$-fractals]{Hilbert-Kunz series, F-signature series, and weak $p$-fractals}


\begin{document}
	\maketitle
	
	\begin{abstract}
		We extend the theory of $p$-fractals of Monsky and Teixeira by introducing the notion of weak $p$-fractal. We prove that for a hypersurface $f$ having rational Hilbert-Kunz series is equivalent to the weak $p$-fractality of the associated function $\phi_{f,p}$ and having rational F-signature series is equivalent to the weak $p$-fractality of the reflection $\overline{\phi}_{f,p}$. In addition, we prove some results characterizing the shape of the generating series of numerical functions which are quasi-polynomials in $p^n$. This is motivated by the fact that the Hilbert-Kunz and F-signature functions take this form in several examples of interest.
	\end{abstract}

	\section{Introduction}

	Let $(R,\mathfrak{m},\KK)$ be a Noetherian local ring of positive characteristic $p$ and assume further that $R$ is a complete domain and the residue field $\KK$ is perfect.
 In this set up, we can consider the following two numerical functions known as \emph{Hilbert-Kunz function} and {\it F-signature function} respectively:
 \[
	\HK_R(n)=\ell_R\left(R/\mathfrak{m}^{[p^n]}\right), \ \ \ \ \ \ \ \ \ \ \FS_R(n)=\frk_R(R^{1/p^n}).
\]
 Here,  $\mathfrak{m}^{[p^n]}$ denotes the Frobenius power of the ideal $\mathfrak{m}$, generated by all the $p^n$-th powers of elements of $\mathfrak{m}$,  $R^{1/p^n}$ denotes the ring obtained by adjoining $p^n$-th roots of elements in $R$, and $\frk_R(R^{1/p^n})$ denotes the maximal rank of a free $R$-summand of $R^{1/p^n}$.
 The Hilbert-Kunz function was first studied by Kunz \cite{Kun69}, while the F-signature function was introduced by Smith and Van den Bergh \cite{SVdB97} in the context of rings with finite F-representation type, and later considered in greater generality by Huneke and Leuschke \cite{HunLeu}, and by Watanabe and Yoshida \cite{WatYosh}.

It is known \cite{Mon83, PolTuc, Tucker12} that both functions grow asymptotically as a polynomial of degree $d=\dim R$ in $p^n$. More precisely, we have
\[
	\HK_R(n)=e_{\HK}(R)p^{dn}+O(p^{(d-1)n}), \ \ \ \ \ \ \ \ \ \ \FS_R(n)=\fs(R)p^{dn}+O(p^{(d-1)n}),
\]
	where $e_{\HK}(R)$ is a positive real number called \emph{Hilbert-Kunz multiplicity}, and $\fs(R)$ is a real number in the interval $[0,1]$ called \emph{F-signature}.
	Like the classical Hilbert-Samuel multiplicity, the Hilbert-Kunz multiplicity and the F-signature encode important information about the singularities of the ring. For example, Watanabe and Yoshida \cite{WatYos2000} proved that, under mild assumptions, $e_{\HK}(R)=1$ if and only if $R$ is a regular local ring, and Huneke and Leuschke \cite{HunLeu} proved that also $\fs(R)=1$ if and only if $R$ is regular. So, the Hilbert-Kunz multiplicity and the F-signature characterize the regularity of the ring in the same way as the Hilbert-Samuel multiplicity does. However, they are not necessarily  integer numbers, and $e_{\HK}(R)$ may be irrational \cite{Bre13, MonConj}.
	
	On the other hand, the Hilbert-Kunz and the F-signature functions may behave very differently from the Hilbert-Samuel function. While the latter is always a polynomial in $n$, the error term $O(p^{(d-1)n})$ of the Hilbert-Kunz and F-signature functions may be very far from being a polynomial in $p^n$. However, if $R$ is  regular in codimension one, then there exists a \emph{second coefficient} for $\HK_R(n)$, that is there exists $\beta\in\mathbb{R}$ such that $\HK_R(n)=e_{\HK}(R)p^{dn}+\beta p^{(d-1)n}+O(p^{(d-2)n})$, and similarly for  $\FS_R(n)$ in the $\mathbb{Q}$-Gorenstein case (see \cite{ChanKurano, HMM04, PolTuc}).
This result is optimal, since there are examples of rings whose Hilbert-Kunz and F-signature functions do not admit a \emph{third coefficient} \cite{CSTZ24, HanMon,MonAlg,MonTrasc,MonTexI,MonTexII}.

Despite the general behavior being quite erratic, for many nice classes of rings the Hilbert-Kunz function has a more regular shape. For example, it is a quasi-polynomial function in $p^n$, i.e. it has the form $a_d(n)p^{dn}+a_{d-1}(n)p^{(d-1)n}+\dots+a_0(n)$ for some periodic rational-valued functions $a_i(n)$, in the following cases: coordinate rings of algebraic curves \cite{Bre07, Tri05}, binomial hypersurfaces \cite{Con96}, affine semigroup rings \cite{Bru05}, ADE singularities \cite{Bri17}. 
For an excellent recent survey on the shape of the Hilbert-Kunz function, see \cite{Chan}.
The F-signature function has been less studied, we mention that it is a quasi-polynomial function in $p^n$ for affine semigroup rings \cite{VKor11} and quotient singularities \cite{CDS19}.
 
In this paper, we concentrate on the  \emph{Hilbert-Kunz series} and the \emph{F-signature series}, that are the generating series 
\[ 
\HKS_R(z)=\sum_{n=0}^{\infty}\HK_R(n)z^n, \ \ \ \ \ \ \ \ \ \  \FSS_R(z)=\sum_{n=0}^{\infty}\FS_R(n)z^n.
\]
Classically, many important algebraic properties of the ring, such as Cohen-Macaulayness, can be read off the  Hilbert-Samuel series which is one of the most important tools that have been employed successfully to study Noetherian local rings. On the other hand, the Hilbert-Kunz and F-signature series are  much more complicated and unpredictable objects. First of all, differently from the Hilbert-Samuel series, they are not always a rational function, i.e., a quotient of polynomials in $\mathbb{Q}[z]$. For the Hilbert-Kunz series, the rationality is known when the ring has finite Cohen-Macaulay type \cite{Sei97}, or when the corresponding function is quasi-polynomial in $p^n$ (Theorem~\ref{thm:succ-serie quasi-pol}). The F-signature series has received even less attention so far being very few the cases where we have a knowledge of the corresponding function.

To study the rationality of the Hilbert-Kunz series of a hypersurface $f\in A=\KK\llbracket x_1,\dots,x_s\rrbracket$,  Monsky and Teixeira \cite{MonTexI, MonTexII} associate to $f$ a function 
\[
 \phi_{f,p}\left(\frac{a}{p^n}\right)=\frac{1}{p^{{sn}}}\dim_{\KK}\big(A/(x_1^{p^n},\dots,x_s^{p^n},f^a)\big).
 \]
Notice that $p^{sn}\phi_{f,p}\left(\frac{1}{p^n}\right)=\HK_{A/f}(n)$ is just the Hilbert-Kunz function of $f$, and the F-signature function can also be recovered from $\phi_{f,p}$ by taking values at points of the form $1-\frac{1}{p^{n}}$ (see Lemma~\ref{lemma_F-signaturefractal}). Moreover $\phi_{f,p}$ is also strictly related to the F-signature of pairs \cite{BST13,CSTZ24}.
Then, Monsky and Teixeira introduce the notion of $p$-fractal by saying that $\phi_{f,p}$ is a $p$-fractal if it lies in a finitely dimensional space stable under the action of certain operators $T_{p^n|b}$ for any $n\in\mathbb{N}$ and integer $0\leq b< p^n$ (see Remark~\ref{rem:pfractalimpliesweakly} for the definition of these operators). They prove that if $\phi_{f,p}$ is a $p$-fractal, then $f$ has rational Hilbert-Kunz series. 

We introduce a weakening of the definition of $p$-fractal, by requiring the stability only under a certain subset of these operators. More precisely, we say that $\phi_{f,p}$ is a \emph{weak $p$-fractal} if  it lies in a finitely dimensional space stable under the operators of the form $T_{p^n|0}$ (Definition~\ref{def:weakp-fractal}). Clearly, $p$-fractal implies weak $p$-fractal, while the converse does not hold (see Example~\ref{ex:weak p-fractal not p-fractal}).
We are able to generalize several properties of $p$-fractals to weak $p$-fractals.
We prove that the set of weak $p$-fractals forms a subalgebra (Lemma~\ref{lem:weak p-frac subalgebra}). This implies that if $f$ and $g$ are power series in disjoint sets of variables such that $\phi_{f,p}$ and $\phi_{g,p}$ are weak $p$-fractals, then $\phi_{fg,p}$ is a weak $p$-fractal (Theorem~\ref{thm:HKofproductfg}).
Moreover, it turns out that the property of being a weak $p$-fractal characterizes the rationality of the corresponding Hilbert-Kunz and F-signature series. More precisely, we are able to prove the following two facts (Corollaries~\ref{cor:HKseriesrational} and \ref{cor:FSS rational iff weak p-fractal}):
\begin{enumerate}
\item $\HKS_{A/f}(z)$ is rational if and only if the function $\phi_{f,p}$ is a weak $p$-fractal.
\item $\FSS_{A/f}(z)$ is rational if and only if the function $\overline{\phi}_{f,p}$ is a weak $p$-fractal, where $\overline{\phi}_{f,p}(t)=\phi_{f,p}(1-t)$ is the reflection of $\phi_{f,p}$ (Definition~\ref{dfn:reflection}).
\end{enumerate}
One may wonder whether the previous two conditions together imply that $\phi_{f,p}$ is a $p$-fractal. However, a conjecture by Monsky suggests that being a $p$-fractal is a stronger condition than having both Hilbert-Kunz and F-signature series rational (Example~\ref{ex:weakpfractalconjecture}).

The remainder of the paper is of a combinatorial nature. Motivated by several examples from Hilbert-Kunz and F-signature functions, we investigate the generating series of numerical functions that are quasi-polynomial in $p^n$. We characterize the form of the generating series by demonstrating that it is always a rational function, where the denominator is a product of factors of the form $z^M-1/p^{jM}$ where $M$ is the common period of the periodic coefficients of the quasi-polynomial (Theorem~\ref{thm:succ-serie quasi-pol}). In the final section of the paper, we focus on more specific quasi-polynomials of the form $a_dp^{nd}+a_0(n)$, where $a_d\in \mathbb{Q}^\times$ and $a_0(n)$ is a $\mathbb{Q}$-valued periodic function.  The main motivation is that in several examples, the Hilbert-Kunz and F-signature functions take this form \cite{Bre07,BreCam18,CDS19,ChanKurano, HMM04,Kreu,Par94, Tri05}. In this case, Theorem~\ref{thm:succ-serie quasi-pol} implies that the denominator of the generating series is a factor of $(1-p^dz)(1-z^M)$. Thus, we inquire whether any simplification between the numerator and the denominator of the series may occur. We show that the factor $1-p^dz$ never cancels out (Lemma~\ref{lemma_1/p^dnoroot}). On the other hand, the numerator of the series may share a root with $1-z^M$, specifically a root of unity. This can also occur in Hilbert-Kunz series, as demonstrated in Example~\ref{ex:simplification of nonprimitive roots} and Example~\ref{ex:rootofunitysimplifies}. Unfortunately, providing a complete characterization of this phenomenon in algebraic or combinatorial terms appears to be complicated. We present some necessary conditions and open questions in the final part of the paper.

\subsection*{Structure of the paper}
In Section~\ref{sec:generatingseries}, we prove some results on the shape of the generating series of numerical functions which are quasi-polynomials in $p^n$. We also provide some elementary but useful results which we rely on in the rest of the paper.
In Section~\ref{section:weakpfractals}, we introduce the notion of weak $p$-fractal, which generalizes $p$-fractals, and we prove that a function $\phi$ is a weak $p$-fractal if and only if the corresponding generating series is rational (Theorem~\ref{thm:weak p-fractal iff series rational}). In Section~\ref{section:HKFsignatureweakpfractals}, we extend some results by Monsky and Teixeira to the more general setting of weak $p$-fractals and prove the rationality of the Hilbert-Kunz and F-signature series of certain hypersurfaces (Theorem~\ref{thm:HKofproductfg}).
In Section~\ref{section:quasipolyna_d+a_0}, we study the shape of the generating series of quasi-polynomial functions in $p^n$ where the leading coefficient is constant and all coefficient functions $a_i$ are zero except the last one, i.e., $e_n=a_dp^{dn}+a_0(n)$.

\subsection*{Acknowledgements}
We would like to thank Alessandro De Stefani, Yusuke Nakajima, and Kevin Tucker for several helpful discussions. We are also grateful to the anonymous referee for a careful reading and for valuable suggestions.

A.~Caminata is supported by the Italian PRIN2022 grants P2022J4HRR ``Mathematical Primitives for Post Quantum Digital Signatures'', and 2022K48YYP ``Unirationality, Hilbert schemes, and singularities'', by the MUR Excellence Department Project awarded to Dipartimento di Matematica, Università di Genova, CUP D33C23001110001, and by the European Union within the program NextGenerationEU.

F. Zerman is supported by the ERC Consolidator grant "Shimura varieties and the BSD conjecture" (grant ID 101001051).


\section{Generating series of quasi-polynomials}\label{sec:generatingseries}
	Let $\{e_n\}_{n\in\mathbb{N}}$ be a sequence of rational numbers. We denote by $$G(e_n;z)=\sum_{n=0}^{\infty}e_nz^n\in \mathbb{Q}\llbracket z\rrbracket$$ the ordinary generating series of $e_n$.
	It is a classical problem to find sufficient and necessary conditions for  $G(e_n;z)$ to be a rational function, i.e., an element of $\mathbb{Q}(z)$. This is the case for example when $e_n$ is a polynomial or quasi-polynomial function in $n$ (see \cite[\S4]{Stan}). In this section, we study the shape of the series $G(e_n;z)$ when $e_n$ is a polynomial or quasi-polynomial function in $p^n$, where $p$ is a fixed prime number. This is motivated by applications to numerical functions in positive characteristic commutative algebra. In fact, it turns out that in many interesting examples the Hilbert-Kunz function and the F-signature function have this shape.

	\begin{definition}
		Let $p$ be a prime number and $d$ be a nonnegative integer. A \emph{quasi-polynomial of degree $d$ in $p^n$} is a sequence $\{e_n\}_{n\in\N}:\mathbb{N}\to\mathbb{Q}$ of the form
		\begin{equation*}
			e_n=a_d(n)p^{dn}+a_{d-1}(n)p^{(d-1)n}+\dots+a_0(n)
		\end{equation*}
		where each $a_i:\mathbb{N}\rightarrow\mathbb{Q}$ is a periodic function with integral period and $a_d(n)$ is not identically zero.
		If  $a_i(n)$ is a constant function for all $i=0,\dots,d$, we say that $e_n$ is a \emph{polynomial of degree $d$ in $p^n$}.
	\end{definition}
	
	We begin with two basic results on generating functions. Recall that a sequence of rational numbers $\{e_n\}_{n\in\mathbb{N}}$ is said to be \emph{linearly recurrent} if there exist $m\in\mathbb{Z}_{>0}$ and $c_j\in\Q$ for $j=1,\dots, m$ such that $e_n=\sum_{j=1}^{m}c_j e_{n-j}$ for all $n$ big enough.
	
	\begin{lemma}	\label{lem:linear recurrence}
		Let $\{e_n\}_{n\in\mathbb{N}}$ be a sequence of rational numbers. Then, $G(e_n;z)\in \Q(z)$ if and only if the sequence $\{e_n\}_{n\in\N}$ is linearly recurrent.
	\end{lemma}
	\begin{proof}
		Suppose that $G(e_n;z)=A(z)/B(z)$ with $A(z), B(z)\in\Q[z]$ and let $B(z)=\sum_{j=0}^{m}b_jz^j$. Then, the relation $B(z)G(e_n;z)=A(z)$ gives the linear recurrence $\sum_{j=0}^{m}b_je_{n-j}=0$ for every $n>\max\{\deg(A(z)), \deg(B(z))\}$. 
  
  Conversely, assume that there is a linear recurrence $\sum_{j=0}^m c_j e_{n-j}=0$ for all sufficiently large $n$ and for fixed $m>0$ and $c_j\in\Q$. Set $B(z)=\sum_{j=0}^m c_jz^j$. Then, for all sufficiently large $n$, the coefficient of $z^n$ inside the product $B(z)G(e_n;z)$ is equal to $\sum_{j=0}^m c_j e_{n-j}=0$, i.e., $B(z)G(e_n;z)$ is a polynomial.
	\end{proof}
	
	\begin{proposition}\label{prop:successione-serie}
		Let $\KK$ be a field, $\{e_n\}_{n\in\N}$ be a sequence of elements of $\KK$ and $\delta_0,\dots,\delta_d\in\KK$ be distinct nonzero elements. The following statements are equivalent:
		\begin{itemize}
			\item[(a)] There exist $a_0,\dots,a_d\in\KK$ such that
			
			\begin{equation*}
				e_n=a_0 \delta_0^n+a_{1}\delta_1^{n}+\dots +a_d\delta_d^n
			\end{equation*}
			for every $n\in\N$.
			\item[(b)] There exist $P(z),Q(z)\in\KK[z]$ such that
			\begin{equation*}
				G(e_n;z)= \frac{P(z)}{Q(z)}
			\end{equation*}
			with $Q(z)=\displaystyle\prod_{i=0}^d(1-\delta_i z)$ and $\deg P(z)<\deg Q(z)=d+1$.
		\end{itemize}
	\end{proposition}
	\begin{proof}
		First, we assume that condition (a) holds. Then
		\begin{equation*}
			\begin{split}
				G(e_n;z)&=\sum_{n=0}^{\infty} e_nz^n=\sum_{n=0}^{\infty} a_0\delta_0^nz^n+\sum_{n=0}^{\infty} a_1\delta_1^nz^n+\dots +\sum_{n=0}^{\infty} a_d\delta_d^nz^n\\
				&=\frac{a_0}{1-\delta_0z}+\frac{a_1}{1-\delta_1z}+\dots+\frac{a_d}{1-\delta_dz},
			\end{split}
		\end{equation*}
		using the formula for the geometric series. Taking the common denominator we easily obtain condition (b). 
		
		Conversely, assume that $e_n$ satisfies (b). Applying the partial fraction decomposition, we have
		\begin{equation*}
			G(e_n;z)=\frac{P(z)}{Q(z)}=\frac{P(z)}{\prod_{i=0}^d(1-\delta_iz)}=\frac{a_0}{1-\delta_0z}+\frac{a_1}{1-\delta_1z}+\dots+\frac{a_d}{1-\delta_dz}
		\end{equation*}
		for some $a_0,\dots,a_d\in\KK$. We define the sequence $	f_n=a_0\delta_0^n+\dots+a_d\delta_d^n$.
		From the previous implication, we see that
		\begin{equation*}
			\sum_{n=0}^{\infty} f_nz^n=\frac{a_0}{1-\delta_0z}+\dots+\frac{a_d}{1-\delta_dz}=\sum_{n=0}^{\infty} e_nz^n
		\end{equation*}
		as elements of $\KK\llbracket z\rrbracket$, hence $e_n=f_n$ for every $n\in\N$.
	\end{proof}
	
	\begin{remark}
		By multiplying $Q(z)$ and $P(z)$ by an appropriate constant, we may write the denominator in point (b) of Proposition~\ref{prop:successione-serie} as $Q(z)=\prod_{i=0}^d(z-1/\delta_i)$.
	\end{remark}

	By applying Proposition \ref{prop:successione-serie} with $\delta_i=p^i$ for every $0\leq i\leq d$ we obtain the following characterization for generating series of polynomial functions in $p^n$.
	
	\begin{corollary}\label{cor:succ-serie pol}
		Let $\{e_n\}_{n\in\N}$ be a sequence of rational numbers and let $d\in\N$. The following statements are equivalent:
		\begin{itemize}
			\item[(a)] There exist $a_0,\dots,a_d\in\Q$ such that $e_n=a_d p^{dn}+a_{d-1}p^{(d-1)n}+\dots +a_0$.
			
			\item[(b)] There exist $P(z),Q(z)\in\Q[z]$ such that
			\begin{equation*}
				G(e_n;z)= \frac{P(z)}{Q(z)}
			\end{equation*}
			with $Q(z)=\prod_{i=0}^d(z-\frac{1}{p^i})$ and $\deg P(z)<\deg Q(z)=d+1$.
		\end{itemize}
		Moreover, when one of the two equivalent conditions is satisfied, then $Q(z)$ can be chosen to be equal to $\prod_i(z-\frac{1}{p^i})$ with $i$ running among the indices such that $a_i\ne 0$.
	\end{corollary}
	
	Now we can state our characterization of quasi-polynomial functions in $p^n$ in terms of generating series. We point out that it is not restrictive to assume that the coefficients $a_i(n)$ of a quasi-polynomial function have the same period $M\in\N_{\ge 1}$. In fact, we may always take $M$ to be the lowest common multiple of the periods of the coefficients.
	
	\begin{theorem}	\label{thm:succ-serie quasi-pol}
		Let $\{e_n\}_{n\in\N}$ be a sequence of rational numbers, let $d\in\mathbb{N}$ and $M\in\N_{\ge 1}$. The following statements are equivalent:
		\begin{itemize}
			\item[(a)] There exist $a_0,\dots,a_d:\mathbb{N}\rightarrow\mathbb{Q}$ periodic functions with period $M\in\N$ such that
			\begin{equation*}
				e_n=a_d(n) p^{dn}+a_{d-1}(n)p^{(d-1)n}+\dots +a_0(n).
			\end{equation*}
			\item[(b)] There exist  $P(z),Q(z)\in\Q[z]$ such that
			\begin{equation*}
				G(e_n;z)= \frac{P(z)}{Q(z)}
			\end{equation*}
			with $Q(z)=\prod_{j=0}^d (z^M-1/p^{jM})$ and $\deg P(z)<\deg Q(z)=M(d+1)$.
		\end{itemize}
		Moreover, let $M_j$ be the minimal period of $a_j(n)$. When one of the two equivalent conditions is satisfied then $Q(z)$ can be chosen to be equal to $\prod_j (z^{M_j}-1/p^{jM_j})$ with $j$ varying among those $j$ such that $a_j(n)$ is not identically zero.
	\end{theorem}
	\begin{proof}
		$(a)\Rightarrow (b)$ Let $e_n$ be as in condition (a). Then, by additivity, we find that
		\begin{equation*}
			G(e_n,z)=\sum_{n=0}^{\infty} e_nz^n=\sum_{n=0}^{\infty} a_d(n)p^{nd}z^n+\dots+\sum_{n=0}^{\infty} a_0(n)z^n.
		\end{equation*}
		Fix $0\le j\le d$ such that $a_j(n)$ is not the zero function and set $f_n^{(j)}=a_j(n)p^{nj}$.  Let $M_j$ be the minimal period of $a_j(n)$ and fix $0\le i <M_j$. The sequence $n\mapsto f_{nM_j+i}^{(j)}$ is of the form
		\begin{equation*}
			\begin{split}
				f_{nM_j+i}^{(j)}=a_j(i)p^{j(nM_j+i)}=a_j(i)p^{ji}(p^{jM_j})^n.
			\end{split}
		\end{equation*}
		Since 
		\begin{equation*}
			G(f_n^{(j)};z)=\sum_{n=0}^\infty f_n^{(j)}z^n=\sum_{i=0}^{M_j-1}\sum_{n=0}^\infty f_{nM_j+i}^{(j)}z^{nM_j+i}=\sum_{i=0}^{M_j-1}\sum_{n=0}^\infty a_j(i)p^{ji}z^i(p^{jM_j}z^{M_j})^n,
		\end{equation*}
		by the properties of the geometric series we conclude that
		\begin{equation}
			\label{equ:serie quasi-polinomi}
			G(f_n^{(j)};z)=\sum_{i=0}^{M_j-1}\frac{a_j(i)p^{ji}z^i}{1-p^{jM_j}z^{M_j}}
		\end{equation}
		Since $G(e_n;z)=\sum_j G(f_n^{(j)};z)$, we have that condition (b) holds. In fact, we also proved that condition (a) implies the "moreover" part of the statement of the proposition. In order to conclude the proof, we are left to show that (b) implies (a).
		
		$(b)\Rightarrow (a)$ Let $\beta_0,\dots,\beta_{d+1}$ be the coefficients of $1,z^M,z^{2M},\dots,z^{(d+1)M}$ in $Q(z)$. We will build a quasi-polynomial $f_n$ such that $\sum f_nz^n=P(z)/Q(z)$, and then conclude that $e_n=f_n$.
		
		Write $P(z)=b_{M(d+1)-1}z^{M(d+1)-1}+\dots+b_1z+b_0$ with $b_j\in\Q$. Since $Q(z)\sum e_nz^n=P(z)$, the numbers $e_0,\dots,e_{(d+1)M-1}$ satisfy the relations
		\begin{equation*}
			\begin{cases}
				&e_i\beta_0=b_i\\
				&e_{M+i}\beta_0+e_i\beta_1=b_{M+i}\\
				&\vdots\\
				&e_{Md+i}\beta_0+e_{M(d-1)+i}\beta_1+\dots+e_i\beta_{d}=b_{Md+i}
			\end{cases}
		\end{equation*}
		for every $0\le i< M$. Let $P_i(z)$ be the unique polynomial of degree $d$ such that $P_i(p^{Mn})=e_{nM+i}$ for every $0\le n\le d$, with $0\le i<M$. We define $f_n=e_n$ if $0\le n\le M(d+1)-1$, and $f_{Mn+i}=P_i(p^{Mn})$ for every $n\ge d+1$.
		It is clear that $f_n$ built in this way is a quasi-polynomial of degree $d$ and period $M$ that coincides with $e_n$ for $n< M(d+1)$. We claim that $f_n=e_n$ for all $n\in\mathbb{N}$.	
		Thanks to the implication $(a)\Rightarrow (b)$, we have that
		\begin{equation*}
			\sum f_nz^n=\frac{R(z)}{Q(z)}
		\end{equation*}
		with $Q(z)$ as in point $(b)$ and $R(z)\in\Q[z]$ such that $\deg R(z)<\deg Q(z)=M(d+1)$. Since $e_n=f_n$ whenever $n< M(d+1)$, we have that
		\begin{equation*}
			Q(z)\sum_{n=0}^{\infty} f_nz^n=R(z)=P(z)=Q(z)\sum_{n=0}^{\infty} e_nz^n.
		\end{equation*}
		This implies that $\sum f_nz^n=\sum e_n z^n$, hence $e_n=f_n$ for every $n\in\N$ as desired.
	\end{proof}
	
	We can think of the power series $G(e_n;z)=\sum e_nz^n$ as a function of complex variable. When $G(e_n;z)=P(z)/Q(z)$ for polynomials $P(z),Q(z)\in\Q[z]$, by the identity principle for holomorphic functions, we have that the function $G(e_n;z)$ is defined everywhere in $\C$ outside the zeroes of $Q(z)$. In particular, when $e_n$ is a polynomial or quasi-polynomial in $p^n$, we know precisely where the zeros of $Q(z)$ may lie thanks to Theorem~\ref{thm:succ-serie quasi-pol} and Corollary~\ref{cor:succ-serie pol}.
	When $e_n$ is not a quasi-polynomial, but it is a \emph{quasi-polynomial up to a certain degree} we can still control the poles of the series $G(e_n;z)$ in a suitable ball centered in the complex origin. We state this precisely in the next proposition.
	
	\begin{proposition}\label{prop:holomorphic quasi-polynomial}
		Let $\{e_n\}_{n\in\N}$ be a sequence of rational numbers.
		Suppose that there exist $d\in\N$ and $0\le \ell\le d$ such that
		\begin{equation*}
			e_n=a_d(n) p^{dn}+a_{d-1}(n)p^{(d-1)n}+\dots +a_\ell(n)p^{\ell n}+f_n
		\end{equation*}
		with $f_n=O(p^{\alpha n})$ for some $\alpha<\ell$ and with $a_j(n)$ periodic of period $M_j\geq1$.
		Then, in the ball $|z|<1/p^{\alpha}$, the holomorphic function $G(e_n;z)$ can only have simple poles and they can only be placed at $M_j$-th roots of $\frac{1}{p^{jM_j}}$, for $\ell\le j\le d$. 
	\end{proposition}
	\begin{proof}
		By Theorem~\ref{thm:succ-serie quasi-pol}
		\begin{equation*}
			\begin{split}
				G(e_n;z)=\sum_{n=0}^\infty e_nz^n=	\sum_{n=0}^\infty (a_d(n) p^{dn}+\dots +a_\ell(n)p^{\ell n})z^n+\sum_{n=0}^\infty f_nz^n=\frac{P(z)}{Q(z)}+\sum_{n=0}^\infty f_nz^n
			\end{split}
		\end{equation*}
		with $Q(z)=\prod_{j=\ell}^d(z^{M_j}-\frac{1}{p^{jM_j}})$ and $P(z)\in\Q[z]$. Since the series $\sum f_nz^n$ is absolutely convergent for $|z|<1/p^{\alpha}$, we have that $Q(z)G(e_n;z)$ is absolutely convergent for $|z|<1/p^{\alpha}$, hence the function $Q(z)G(e_n;z)$ is holomorphic for $|z|<1/p^{\alpha}$. This implies that, inside the ball $|z|<1/p^{\alpha}$, the function $G(e_n;z)$ can just have poles in the roots of $Q(z)$, and they can be at most simple poles. 
	\end{proof}
	
	In the following corollary, we specialize the result of the previous proposition to the case when the leading part of $e_n$ has a polynomial shape. This happens for the Hilbert-Kunz and the F-signature function which are known to be of the form $\mu p^{nd}+O(p^{n(d-1)})$ where $d$ is the dimension of the ring and $\mu$ is the Hilbert-Kunz multiplicity (resp. the F-signature).

	\begin{corollary}
		\label{cor:poles of gen funct in a ball}
		Let $\{e_n\}_{n\in\N}$ be a sequence of rational numbers.
		Suppose that there exist $d\in\N$, $a_0,\dots,a_d\in\Q$ nonzero and $0\le \ell\le d$ such that 
		\begin{equation*}
			e_n=a_d p^{dn}+a_{d-1}p^{(d-1)n}+\dots +a_\ell p^{\ell n}+f_n
		\end{equation*}
		with $|f_n|=O(p^{\alpha n})$ for some $\alpha<\ell$.
		Then, in the ball $|z|<1/p^{\alpha}$, the holomorphic function $G(e_n;z)$ can only have simple poles and they can only be places at $z=\frac{1}{p^{i}}$ for every $\ell\le i\le d$. 
		Moreover, we have that
		\begin{equation*}
			-a_ip^{-i}=\mathop{\mathrm{Res}}_{z=1/p^{i}}G(e_n;z)=\lim_{z\to 1/p^i}\bigg(z-\frac{1}{p^i}\bigg)G(e_n;z)
		\end{equation*}
		for every $\ell\le i\le d$.
	\end{corollary}
	\begin{proof}
		The first part is Proposition \ref{prop:holomorphic quasi-polynomial} when $M_j=1$ for all $j$.
		Then, as seen in the proof of Proposition \ref{prop:holomorphic quasi-polynomial}, the function $G(e_n;z)$ can be written as
		\begin{equation*}
			\frac{-a_dp^{-d}}{z-1/p^{d}}+\dots+\frac{-a_\ell p^{-\ell}}{z-1/p^{\ell}}+g(z)
		\end{equation*}
		whenever $|z|<1/p^{\alpha}$, where $g(z)$ is a holomorphic function in this ball. The claim follows.
	\end{proof}
	
	\begin{remark}\label{rem:Hkmultfromseries}
		If $e_n$ is the Hilbert-Kunz function of a $d$-dimensional local ring $R$, then the previous corollary tells us that we can compute the Hilbert-Kunz multiplicity $\eHK(R)$ from the generating series $\HKS_R(z)$ by 
		\[
		\eHK(R)=-{p^d}\cdot\mathop{\mathrm{Res}}_{z=1/p^{d}}\HKS_R(z)=\lim_{z\to 1/p^i}(1-p^dz)\HKS_R(z). 
		\]
		In particular, when $\HKS_R(z)\in\Q(z)$ is a rational function, then $\eHK(R)\in\Q$ as observed also by Seibert \cite[Corollary~2.7]{Sei97}. Similarly we can compute the F-signature $\fs(R)$ from the generating series $\FSS_R(z)$ by
        \begin{equation*}
            \fs(R)=-{p^d}\cdot\mathop{\mathrm{Res}}_{z=1/p^{d}}\FSS_R(z)=\lim_{z\to 1/p^i}(1-p^dz)\FSS_R(z). 
        \end{equation*}
        and also conclude that $\fs(R)\in \Q$ when  $\FSS_R(z)\in\Q(z)$.
	\end{remark}


	\section{Weak $p$-fractals}\label{section:weakpfractals}
In order to study the rationality of the Hilbert-Kunz series of a hypersurface $f$ in characteristic $p>0$, Monsky and Teixeira \cite{MonTexI, MonTexII} associate to it a function $\phi_{f,p}$, they develop the theory of \emph{$p$-fractals} and prove that if $\phi_{f,p}$ is a $p$-fractal then the Hilbert-Kunz series for $f$ is rational. In this section, we are going to generalize the theory of $p$-fractals in order to characterize the rationality of the Hilbert-Kunz and F-signature series in terms of specific properties of  \enquote{fractality} of some $\phi_{f,p}$-type functions.

We begin with some notation. 

\begin{definition}\label{def:weakp-fractal}
	Let $\mathscr{V}=\left\{\frac{1}{p^m}: \ m\in\N\right\}$ and let $\Q^{\mathscr{V}}$ be the space of functions from $\mathscr{V}$ to $\Q$. 	An element $\phi$ of $\Q^{\mathscr{V}}$ is called a \emph{weak $p$-fractal} if $\phi$ lies in a finitely dimensional subspace of $\Q^{\mathscr{V}}$ stable under the set of operators $\{T_{p^n|0} \ | \ n\in\N\}$, where
	\begin{equation*}
		(T_{p^n|0}\phi)\left(\frac{1}{p^m}\right)=\phi\left(\frac{1}{p^{m+n}}\right)
	\end{equation*}
 for any $\phi\in \Q^{\mathscr{V}}$ and $m\in\N$.
\end{definition}

Notice that in order to show that $\phi\in \Q^{\mathscr{V}}$ is a weak $p$-fractal it is enough to show that $\phi$ lies in a finitely generated subspace of $\Q^{\mathscr{V}}$ stable under the operator $T_{p|0}$.

\begin{remark}\label{rem:pfractalimpliesweakly}
	Let now $\mathscr{I}=[0,1]\cap\left\{\frac{a}{p^m}: \ a,m\in\N\right\}$. In \cite{MonTexI,MonTexII} Monsky and Teixeira define operators $T_{p^n|b}$ for any $n\in\mathbb{N}$ and integer $0\leq b<p^n$ on $\Q^{\mathscr{I}}$ as 
	\[
	(T_{p^n|b}\phi)\left(\frac{a}{p^m}\right)=\phi\left(\frac{a+bp^m}{p^{m+n}}\right).
	\]
	They call a function $\phi\in \Q^{\mathscr{I}}$ a \emph{$p$-fractal} if the set of all $T_{p^n|b}\phi$ with $n\in\mathbb{N}$ and  $0\leq b<p^n$ span a finitely dimensional $\Q$-subspace of $\Q^{\mathscr{I}}$ .
	Since we have an inclusion $\mathscr{V}\subseteq\mathscr{I}$, we can restrict any function $\phi\in \Q^{\mathscr{I}}$ to a function in $\Q^{\mathscr{V}}$, which we still denote by $\phi$. Similarly, the restriction of the operators $T_{p^n|0}$ defined here coincide with the operators of Definition~\ref{def:weakp-fractal}.
	Then, it is clear that 
	\begin{center}
		$\phi$ $p$-fractal $\Longrightarrow$ $\phi$ weak $p$-fractal.
	\end{center}
	The converse implication does not hold as we will see in Example~\ref{ex:weak p-fractal not p-fractal}.
\end{remark}

We can endow the space of functions $\mathbb{Q}^\mathscr{V}$ with an algebra structure by pointwise sum and multiplication. It turns out that the set of weak $p$-fractals forms a subalgebra of $\Q^{\mathscr{V}}$. This extends the analogous result for $p$-fractals (see \cite[Remark 2.2]{MonTexI}). 

\begin{lemma}
	\label{lem:weak p-frac subalgebra}
	The set of weak $p$-fractals forms a subalgebra of $\Q^{\mathscr{V}}$.
\end{lemma}

\begin{proof}
	Call $X=\Q^{\mathscr{V}}$. Let $\phi,\psi\in X$ be weak $p$-fractals and let $V_\phi$ and $V_\psi$ be two finitely dimensional subspaces of $X$ containing respectively $\phi$ and $\psi$, both of them stable under the operators $T_{p^n|0}$. The linearity of $T_{p^n|0}$ implies that also $V_\phi+ V_\psi$ is stable under $T_{p^n|0}$, hence $\lambda\phi+\mu\psi$ is a weak $p$-fractal for any $\lambda,\mu\in\Q$. Therefore, the set of weak $p$-fractals is a linear subspace of $X$.
	
	Let $V_\phi\cdot V_\psi$ be the linear subspace of $X$ spanned by the elements $\alpha\cdot\beta$ with $\alpha\in V_\phi$ and $\beta\in V_\psi$. This space is finitely dimensional, since it can be spanned by all the combinations of products of an element of a basis of $V_\phi$ with an element of a basis of $V_\psi$. Moreover, since $T_{p^n|0}$ is linear with respect to sum and product, we have that $V_\phi\cdot V_\psi$ is stable under the action of $T_{p^n|0}$. Since $\phi\psi\in V_\phi\cdot V_\psi$, we conclude that $\phi\psi$ is a weak $p$-fractal.
\end{proof}

\begin{definition}\label{dfn:sequence attached to phi}
	Given  an element $\phi\in \Q^{\mathscr{V}}$ and a non-negative integer $s$, we can associate to $\phi$ a sequence of rational numbers $\{e_{s,n}(\phi)\}_{n\in\N}$ as follows
	\[
	e_{s,n}(\phi)= p^{ns}(T_{p^n|0}\phi)(1)=p^{ns}\phi\left(\frac{1}{p^n}\right)\in\Q.
	\]
	When the integer $s$ will be clear from the context, we will just write $e_{n}(\phi)=	e_{s,n}(\phi)$.
\end{definition}

We can characterize the rationality of the generating series $G(e_{s,n}(\phi);z)$ in terms of properties of $\phi$. This will be useful in the next sections since sequences $e_{s,n}(\phi)$ associated to certain functions $\phi$ coincide with the Hilbert-Kunz and F-signature functions of hypersurfaces.
The following characterization of rationality generalizes \cite[Proposition~1.2]{MonTexII}.

\begin{theorem}
	\label{thm:weak p-fractal iff series rational}
	Let $\phi\in \Q^{\mathscr{V}}$ and let $s\ge 0$ be an integer. Then $\phi$ is a weak $p$-fractal if and only if $G(e_{s,n}(\phi);z)\in\Q(z)$.
\end{theorem}
\begin{proof}
		Let $T=T_{p|0}$ and $X=\Q^{\mathscr{V}}$. If $\phi$ is a weak $p$-fractal, the $\Q$-subspace $V$ of $X$ spanned by $T^n(\phi)$ for $n=0,1,2,\dots$ is finite dimensional, and $T$ maps $V\to V$. By Cayley–Hamilton theorem, the restriction of $T$ to $V$ satisfies a polynomial identity of the form $T^\ell=c_1T^{\ell-1}+c_2T^{\ell-2}+\dots+c_\ell$ for some $c_i\in\Q$ and $\ell\in\N$. This equation can be applied to any $T^n(\phi)$ and shows that the sequence $(T^n(\phi))_{n\in\N}$ is linearly recurrent. Indeed 
	\begin{equation*}
		T^{\ell+n}(\phi)=c_1T^{\ell+n-1}(\phi)+c_2 T^{\ell+n-2}(\phi)+\dots+c_\ell T^{n}(\phi)
	\end{equation*}
	for any $n\in\N$.  By multiplying both sides by $p^{(\ell+n)s}$ we obtain the recurrence
	\begin{equation*}
		p^{(\ell+n)s}T^{\ell+n}(\phi)=c_1p^{s}p^{(\ell+n-1)s}T^{\ell+n-1}(\phi)+\dots+c_\ell p^{\ell s}p^{ns}T^n(\phi)
	\end{equation*}
	for any $n\in\N$.
	Evaluating at $1$, we get a linear recursion for the sequence $\{p^{ns}T^n(\phi)(1)\}_{n\in\N}=\{e_{s,n}(\phi)\}_{n\in\N}$ for any $n\ge \ell$. By applying Lemma~\ref{lem:linear recurrence}, we conclude that $G(e_{s,n}(\phi);z)$ is a rational function.
	
	Conversely, assume that $G(e_{s,n}(\phi)(1);z)$ is a rational function. Again by Lemma~\ref{lem:linear recurrence} the sequence $\{p^{ns}T^n(\phi)(1)\}_{n\in \N}$ is linearly recurrent, i.e., there exist $k, \ell\in\N$ and $c_j\in\Q$ for $j=0,1,\dots,\ell$ such that
	\begin{equation*}
		p^{s(\ell+n)}T^{\ell+n}(\phi)(1)=c_1p^{s(\ell+n-1)}T^{\ell+n-1}(\phi)(1)+c_2 p^{s(\ell+n-2)}T^{\ell+n-2}(\phi)(1)+\dots+c_\ell p^{sn}T^{n}(\phi)(1)
	\end{equation*}
	for any $n\ge k$. Dividing both sides by $p^{s(\ell+n)}$ we find the linear recurrence
	\begin{equation*}
		T^{\ell+n}(\phi)(1)=d_1T^{\ell+n-1}(\phi)(1)+d_2T^{\ell+n-2}(\phi)(1)+\dots+d_\ell T^{n}(\phi)(1)
	\end{equation*} 
	for any $n\ge k$ and some appropriate (but fixed) $d_1,\dots,d_\ell\in\Q$. Since $T^{n-k}(\phi)(1)=\phi(1/p^{n-k})$, the above relation is equivalent to 
	\begin{equation*}
		T^{\ell+k}(\phi)\Big(\frac{1}{p^{n-k}}\Big)=d_1T^{\ell+k-1}(\phi)\Big(\frac{1}{p^{n-k}}\Big)+d_2 T^{\ell+k-2}(\phi)\Big(\frac{1}{p^{n-k}}\Big)+\dots+d_\ell T^k(\phi)\Big(\frac{1}{p^{n-k}}\Big)
	\end{equation*}
	for every $n\ge k$.
	This means that, for every $t\in\mathscr{V}$, we have that
	\begin{equation*}
		T^{\ell+k}(\phi)(t)=d_1T^{\ell+k-1}(\phi)(t)+d_2 T^{\ell+k-2}(\phi)(t)+\dots+d_\ell T^k(\phi)(t).
	\end{equation*}
	This implies that the $\Q$-subspace $V$ of $\Q^{\mathscr{V}}$ generated by $\{T^n(\phi)$ for $n\in\N\}$ is in fact generated by $\{T^n(\phi)$ for $n<\ell+k\}$, hence its dimension is finite. Since $V$ is $T$-stable, we conclude that $\phi$ is a weak $p$-fractal.
\end{proof}


Keep the notations of Remark \ref{rem:pfractalimpliesweakly}. After \cite{MonTexI,MonTexII}, we introduce the following definition.

\begin{definition}
\label{dfn:reflection}
	For any $\phi\in\Q^{\mathscr{I}}$ we define the \emph{reflection} $\overline{\phi}$ of $\phi$ as
	\begin{equation*}
		\overline{\phi}(t)=\phi(1-t)
	\end{equation*}
	for every $t\in\mathscr{I}$.
\end{definition}

\begin{remark}\label{rem:pfractalreflection}
	Let $\phi\in\Q^\mathscr{I}$. If $\phi$ is a $p$-fractal, then $\overline{\phi}$ is a weak $p$-fractal. This is because, in general, $\phi$ is a $p$-fractal if and only if $\overline{\phi}$ is a $p$-fractal. This can be checked directly by observing that $T_{p^n|a}\overline{\phi}=\overline{T_{p^n|p^n-a-1}\phi}$.
\end{remark}

The reflection operator will play a fundamental role in the next sections. It also gives an easy way to find examples of weak $p$-fractals that are not $p$-fractals.

\begin{example}\label{ex:weak p-fractal not p-fractal}
    Let $\phi\in\Q^{\mathscr{I}}$ be the function defined as
    \begin{equation*}
        \phi\bigg(\frac{a}{p^n}\bigg)=\begin{cases}
            1&\text{if $a=p^n-1$ and $n$ is a perfect square}\\
            0&\text{otherwise}.
        \end{cases}
    \end{equation*}
    Since $e_{0,n}(\phi)=0$ for every $n\ge 2$, Theorem \ref{thm:weak p-fractal iff series rational} implies that $\phi$ is a weak $p$-fractal. On the other hand, we have that
     \begin{equation*}
        e_{0,n}(\overline{\phi})=\phi\bigg(\frac{p^n-1}{p^n}\bigg)=\begin{cases}
            1&\text{if $n$ is a perfect square}\\
            0&\text{otherwise}.
        \end{cases}
    \end{equation*}
    The sequence $e_{0,n}(\overline{\phi})$ is not linearly recurrent, since the difference between two adjacent perfect squares goes to infinity when they approach infinity. Therefore, by Lemma \ref{lem:linear recurrence}, the series $G(e_{0,n}(\overline{\phi});z)$ is not rational and Theorem \ref{thm:weak p-fractal iff series rational} yields that $\overline{\phi}$ is not a weak $p$-fractal. Then, the considerations of Remark \ref{rem:pfractalreflection} imply that $\phi$ is not a $p$-fractal.
\end{example}

\section{Hilbert-Kunz, F-signature, and weak $p$-fractals}\label{section:HKFsignatureweakpfractals}
 
	Let $\KK$ be a perfect field of characteristic $p>0$, let $A=\KK\llbracket x_1,\dots,x_s\rrbracket$ be the power series ring in $s$ variables and consider an element $f$ of the maximal ideal of $A$. We are interested in the Hilbert-Kunz and F-signature functions and series of the quotient $A/f$. 
	\begin{notation}
	From now on, we will denote $\HK_{A/f}(n)$ simply by $\HK_f(n)$ and similarly we will use $\HKS_f(z)$, $e_{\HK}(f)$, $\FS_f(n)$, $\FSS_f(z)$ and $\fs(f)$ in place of writing explicitly the quotient $A/f$.
	\end{notation}
	
After \cite{MonTexI,MonTexII}, we associate to any element $f\in A$ an element  $\phi_{f,p}$ of the space $\Q^{\mathscr{I}}$ 	of functions from $\mathscr{I}=[0,1]\cap\left\{\frac{a}{p^m}: \ a,m\in\N\right\}$ to $\Q$ 
	
\begin{equation*}
	\begin{split}
			\phi_{f,p}:
			\frac{a}{p^n}\longmapsto p^{-sn}\dim_{\KK}\left(A/(x_1^{p^n},\dots,x_s^{p^n},f^a)\right).
		\end{split}
\end{equation*}
As observed in Remark~\ref{rem:pfractalimpliesweakly}, the restriction of $\phi_{f,p}$ to $\mathscr{V}=[0,1]\cap\left\{\frac{1}{p^m}: \ m\in\N\right\}\subseteq\mathscr{I}$ gives an element of the space $\Q^{\mathscr{V}}$ which we still denote by $\phi_{f,p}$.

\begin{remark}\label{rem:HK=phi_f}
	The sequence of Definition \ref{dfn:sequence attached to phi}  associated with the function $\phi_{f,p}$ is precisely the Hilbert-Kunz function of the quotient ring $A/f$ with respect to the maximal ideal $(x_1,\dots,x_s)$,
	that is
	\[
\HK_{f}(n)=	e_{n}(\phi_{f,p}) \ \ \forall n\in\N.
	\]			
	\end{remark}
	
	The previous remark immediately yields the following corollary of Theorem~\ref{thm:weak p-fractal iff series rational}.
	
	\begin{corollary}\label{cor:HKseriesrational}
		The Hilbert-Kunz series of $A/f$ is rational if and only if the function $\phi_{f,p}$ is a weak $p$-fractal.		
	\end{corollary}

The reflection (see Definition \ref{dfn:reflection}) of the function $\phi_{f,p}$ is connected with the F-signature function of $A/f$. The proof of the following result can be found also in \cite[Proposition~4.1]{BST13}

\begin{lemma}\label{lemma_F-signaturefractal}
	For any $f\in A$ and any $n\in\N$, we have
	\[
	\FS_{f}(n)=p^{sn}-e_{n}(\overline{\phi}_{f,p}).
	\]	
\end{lemma}

\begin{proof}
By \cite{Tucker12}, we can compute the F-signature function of $R=A/f$ as $\FS_{f}(n)=\dim_\KK\left(R/I_n\right)$, where $I_n$ is the ideal of $R$ defined as
\[
I_n=\{r\in R: \varphi(r^{1/p^n})\in\mathfrak{m} \ \forall\varphi\in\mathrm{Hom}_R(R^{1/p^n},R)\},
\]
with $\mathfrak{m}=(x_1,\dots,x_s)$ the maximal ideal of $R$.
We can write the ideal $I_n$ also as
\[
I_n=\frac{\mathfrak{m}^{[p^n]}:\left((f)^{[p^n]}:f\right)}{(f)}=\left(\mathfrak{m}^{[p^n]}:f^{p^n-1}\right)/(f),
\]
where $\mathfrak{m}^{[p^n]}=(x_1^{p^n},\dots,x_s^{p^n})$ is the Frobenius power of $\mathfrak{m}$.
Moreover, multiplication by $f^{p^n-1}$ induces a short exact sequence
\[
0\rightarrow A/\mathfrak{m}^{[p^n]}:f^{p^{n}-1}\longrightarrow A/\mathfrak{m}^{[p^n]}\longrightarrow A/(\mathfrak{m}^{[p^n]},f^{p^n-1})\rightarrow0.
\]
Since the $\KK$-vector space dimension of the quotient of $R$ modulo the ideal $\left(\mathfrak{m}^{[p^n]}:f^{p^n-1}\right)/(f)$ is the same as the $\KK$-vector space dimension of the quotient of $A$ modulo the ideal $\left(\mathfrak{m}^{[p^n]}:f^{p^n-1}\right)$,  we obtain
\[
\begin{split}
	\FS_{f}(n)&=	\dim_\KK \left(A/\left(\mathfrak{m}^{[p^n]}:f^{p^n-1}\right)\right)
	= \dim_\KK \left(A/(\mathfrak{m}^{p^n})\right)-\dim_\KK\left(A/\left(\mathfrak{m}^{[p^n]},f^{p^n-1}\right)\right)\\
	&= p^{sn}-p^{sn}\phi_{f,p}\left(\frac{p^n-1}{p^n}\right)
	=p^{sn}-p^{sn}\phi_{f,p}\left(1-\frac{1}{p^n}\right)\\
	&=p^{sn}-p^{sn}\overline{\phi}_{f,p}\left(\frac{1}{p^n}\right)
	=p^{sn}-e_{n}(\overline{\phi}_{f,p})
\end{split}
\]
as required.
\end{proof}

The previous lemma together with Theorem~\ref{thm:weak p-fractal iff series rational} gives the following result about the rationality of the F-signature series. As usual, we denote the restriction of a function $\overline{\phi}_{f,p}\in \Q^\mathscr{I}$ to $\mathscr{V}$ again by $\overline{\phi}_{f,p}$.

\begin{corollary}
		\label{cor:FSS rational iff weak p-fractal}
		The F-signature series of $A/f$ is rational if and only if the function $\overline{\phi}_{f,p}$ is a weak $p$-fractal.
\end{corollary}

\begin{proof}
By Lemma~\ref{lemma_F-signaturefractal} we have
\[
\FSS_{f}(z)=\sum_{n=0}^{\infty}\FS_{f}(n)z^n=\sum_{n=0}^{\infty}\left(p^{sn}-e_{n}(\overline{\phi}_{f,p})\right)z^n =\sum_{n=0}^{\infty}p^{sn}z^n-\sum_{n=0}^{\infty}e_{n}(\overline{\phi}_{f,p})z^n.
\]
Since $\sum_{n=0}^{\infty}p^{sn}z^n\in\Q(z)$, the claim follows from Theorem~\ref{thm:weak p-fractal iff series rational}.
\end{proof}

	A consequence of Corollaries~\ref{cor:HKseriesrational} and \ref{cor:FSS rational iff weak p-fractal} is that a hypersurface $A/f$ has both Hilbert-Kunz and F-signature series rational if and only if both $\phi_{f,p}$ and $\overline{\phi}_{f,p}$ are weak $p$-fractals. In particular, if $\phi_{f,p}$ is a $p$-fractal then by Remark~\ref{rem:pfractalreflection} both the Hilbert-Kunz series and the F-signature series of $A/f$ are rational.
	However, $\phi_{f,p}$ being a $p$-fractal seems a stronger condition than having both Hilbert-Kunz and F-signature series rational. Here we give a conjectural example.

\begin{example}\label{ex:weakpfractalconjecture}
		Let $p=2$, the element $f=x^3+y^3+xyz\in A=\mathbb{F}_2\llbracket x,y,z\rrbracket$ has rational F-signature series and rational Hilbert-Kunz series (i.e. $\phi_{f,p}$ and $\overline{\phi}_{f,p}$ are weak $p$-fractals), but $\phi_{f,p}$ is (conjecturally) not a $p$-fractal. 
		The fact that $\HKS_{f}(z)$ is rational follows from the fact that it defines a nodal cubic, hence its Hilbert-Kunz function is a quasi-polynomial in $p^n$ (see \cite[Theorem 3]{BuchC97}). 
		We show that the F-signature function of $A/f$ is $\FS_{f}(n)=1$ for all $n\in\mathbb{Z}_+$, so $\FSS_{f}(z)=\frac{1}{1-z}\in\Q(z)$.
		By the proof of Lemma~\ref{lemma_F-signaturefractal}, it is enough to compute the length of the quotient $A/(x^{p^n},y^{p^n},z^{p^n},f^{p^n-1})$.
		First, we claim that 
		\[
		A/\big(x^{p^n},y^{p^n},z^{p^n},f^{p^n-1}\big) = A/\big(x^{p^n},y^{p^n},z^{p^n},(xyz)^{p^n-1}\big).
		\]
		In fact, since $f$ is homogeneous of degree $3$, the polynomial $f^{p^n-1}$ is homogeneous of degree $3(p^n-1)$. The only monomial of this degree that does not lie in $\mathfrak{m}^{[p^n]}=(x^{p^n},y^{p^n},z^{p^n})$ is $(xyz)^{p^n-1}$, which appears with coefficient 1 in $f^{p^n-1}$. Therefore, it does not vanish modulo $p$ and the claim is proved. Finally, by  Lemma~\ref{lemma_F-signaturefractal} we have
		\[
		\FS_{f}(n)=p^{3n}-\dim_\KK \left( A/\big(x^{p^n},y^{p^n},z^{p^n},(xyz)^{p^n-1}\big)\right)=p^{3n}-(p^{3n}-1)=1.
		\]
		Now, assuming Monsky's conjecture \cite{MonConj}, \cite[Conjecture 2.3]{MonAlg}, we show that the function $\phi_{f,p}$ is not a $p$-fractal.  Assume by contradiction that $\phi_{f,p}$ is a $p$-fractal. From \cite[Theorem 1]{MonTexI} and \cite[Proposition 4.3]{MonTexII} we obtain that the element  $f+uv\in\mathbb{F}_2\llbracket x,y,z,u,v\rrbracket$ has rational Hilbert-Kunz series, hence its Hilbert-Kunz multiplicity is rational by Remark~\ref{rem:Hkmultfromseries}. On the other hand, Monsky, assuming the above mentioned conjecture, proves that the Hilbert-Kunz multiplicity of $f+uv$ is not rational (see the lines after \cite[Theorem 2.4]{MonAlg}).
	\end{example}

Given two power series $f\in\KK\llbracket x_1,\dots,x_s\rrbracket$ and $g\in\KK\llbracket y_1,\dots,y_{s'}\rrbracket$, we consider the power series $f\cdot g$ in the power series ring in $s+s'$ variables $\KK\llbracket x_1,\dots,x_s,y_1,\dots,y_{s'}\rrbracket$. We can extend some results from \cite[\S4]{MonTexII} about the $p$-fractality of the function $\phi_{f g,p}$ to the setting of weak $p$-fractals. The key observation is the following formula
\begin{equation}\label{eq:phi_fg}
	\phi_{fg,p}=\phi_{f,p}+\phi_{g,p}-\phi_{f,p}\phi_{g,p}
\end{equation}
that can be found in the proof of \cite[Proposition~4.2]{MonTexII}.

\begin{theorem}\label{thm:HKofproductfg}
Let $f\in\KK\llbracket x_1,\dots,x_s\rrbracket$ and $g\in\KK\llbracket y_1,\dots,y_{s'}\rrbracket$ be two power series and let $fg\in\KK\llbracket x_1,\dots,x_s,y_1,\dots,y_{s'}\rrbracket$ be their product. Then, the following facts hold.
\begin{enumerate}[(i)]
\item If $\phi_{f,p}$ and $\phi_{g,p}$ are weak $p$-fractals, then  $\phi_{fg,p}$ is a weak $p$-fractal.
\item If $f$ and $g$ have rational Hilbert-Kunz series (resp., F-signature series), then  $f g$ has rational Hilbert-Kunz series (resp., F-signature series).
\item If $f$ and $g$ have quasi-polynomial Hilbert-Kunz function (resp., F-signature function), then $f g$ has quasi-polynomial Hilbert-Kunz function (resp., F-signature function).
\item $e_{\HK}(fg)=e_{\HK}(f)+e_{\HK}(g)$.
\item $\fs(fg)=0$.
\end{enumerate}
\end{theorem}

\begin{proof}
	\begin{enumerate}[(i)]
		\item It follows from the formula \eqref{eq:phi_fg} and the fact that weak $p$-fractals form an algebra by Lemma~\ref{lem:weak p-frac subalgebra}.
		\item The statement about Hilbert-Kunz series follows immediately from \textrm{(i)} and Corollary~\ref{cor:HKseriesrational}. Turning to F-signature series, from \eqref{eq:phi_fg} and the definition of reflection, we know that
		\begin{equation}\label{equ:reflection of phi_fg}
			\overline{\phi}_{fg,p}=\overline{\phi_{f,p}+\phi_{g,p}-\phi_{f,p}\phi_{g,p}}=\overline{\phi}_{f,p}+\overline{\phi}_{g,p}-\overline{\phi}_{f,p}\overline{\phi}_{g,p}.
		\end{equation}
		Now, $\overline{\phi}_{f,p}$ and $\overline{\phi}_{g,p}$ are weak $p$-fractal by assumption, so $\overline{\phi}_{fg,p}$ is also a weak $p$-fractal by Lemma~\ref{lem:weak p-frac subalgebra}. So the F-signature series of $fg$ is rational by Corollary~\ref{cor:FSS rational iff weak p-fractal}.
		\item The statement about Hilbert-Kunz functions (resp., F-signature functions) is a consequence of \eqref{eq:phi_fg} and Remark \ref{rem:HK=phi_f} (resp., \eqref{equ:reflection of phi_fg} and Lemma \ref{lemma_F-signaturefractal}), and the fact that sums and products of periodic functions are again periodic functions. 
		\item We know that 
		\begin{equation*}
			\HK_f(n)=e_{\HK}(f)p^{n(s-1)}+\text{terms of degree less then $s-1$ in $p^n$},
		\end{equation*}
		and similarly for $\HK_g(n)$. This is equivalent to saying that
		\begin{equation*}
			\phi_{f,p}(t)=e_{\HK}(f)t+\text{terms of degree greater than $1$ in $t$}
		\end{equation*}
		for $t\in\mathscr{I}$, and similarly for $\phi_g$. Using the formula \eqref{eq:phi_fg}, we see that
		\begin{equation*}
			\phi_{fg,p}(t)=(e_{\HK}(f)+e_{\HK}(g))t+\text{terms of degree greater than $1$ in $t$}
		\end{equation*}
		for $t\in\mathscr{I}$. This is equivalent to saying that $e_{\HK}(fg)=e_{\HK}(f)+e_{\HK}(g)$.
            \item This can be proved in the same way as \textrm{(iv)}, using Lemma \ref{lemma_F-signaturefractal} and $\overline{\phi}$ instead of $\phi$. Alternatively, one can notice that the quotient ring $\KK\llbracket x_1,\dots,y_{s'}\rrbracket/(fg)$ is not a domain, thus not strongly F-regular. Hence its F-signature is zero.
	\end{enumerate}
\end{proof}


\section{Quasi-polynomials of the form $a_dp^{dn}+a_0(n)$}\label{section:quasipolyna_d+a_0}
	
	In this section, we study more in details the behavior of the generating series $G(e_n;z)$ for sequences of type $e_n=a_dp^{dn}+a_0(n)$ where $d\ge 1$, $a_d\in\Q^\times$ and $a_0(n)$ is a $\Q$-valued periodic function in $n\in\N$ with period $M\geq1$. This is motivated by the fact that in several examples the Hilbert-Kunz function takes this shape. This is the case for two large classes of rings: one-dimensional standard-graded rings \cite{Kreu} and two-dimensional normal, F-finite domains \cite{ChanKurano, HMM04}. The latter includes for example homogeneous coordinate rings of smooth projective curves \cite{Bre07, BreCam18, Tri05}. The Hilbert-Kunz function of some other special families of rings have this shape, such as cuspidal plane cubics \cite{Par94}. Moreover, the F-signature function of certain invariant rings also have this shape \cite{CDS19}.
 
 \subsection{The shape of $G(e_n;z)=P(z)/Q(z)$ }
	From Theorem~\ref{thm:succ-serie quasi-pol} and Equation (\ref{equ:serie quasi-polinomi}), we have
	\begin{equation*}
		\begin{split}
			G(e_n;z)&=\frac{a_d}{1-p^d z}+\sum_{i=0}^{M-1}\frac{a_0(i)z^i}{1-z^M}=\\
			&=\frac{(-a_d-p^da_0(M-1))z^M+\sum_{i=1}^{M-1}(a_0(i)-p^da_0(i-1))z^i+(a_d+a_0(0))}{(1-p^dz)(1-z^M)}.
		\end{split}
	\end{equation*}
	We call $P(z)$ the numerator and $Q(z)$ the denominator of this fraction. Our goal is to study when it is possible that $P(z)$ and $Q(z)$ have common factors, so that some cancellation in the rational function  $G(e_n;z)=P(z)/Q(z)$  occurs.
	
	\begin{lemma}\label{lemma_1/p^dnoroot}
		The number $1/p^d$ is not a root of $P(z)$.
	\end{lemma}
	\begin{proof}
		Assume by contradiction that  $P(1/p^d)=0$, then
		\begin{equation*}
			\begin{split}
				0&=(-a_d-p^da_0(M-1))p^{-dM}+\sum_{i=1}^{M-1}(a_0(i)-p^da_0(i-1))p^{-di}+(a_d+a_0(0))=\\
				&=-a_dp^{-dM}+a_d\ne 0,
			\end{split}
		\end{equation*}
		which is impossible.
	\end{proof}
	
	We examine when it is possible that $P(z)$ is not coprime with $(1-z^M)$. Observe that if $\zeta$ is an $M$-th root of 1, then $P(\zeta)$ only depends on the values of $a_0$, because the terms containing $a_d$ simplify. This means that the condition for $\zeta$ of being a root of $P(z)$ only depends on the values of the periodic function $a_0$.
	
	A straightforward computation shows that $P(1)=0$ if and only if $\sum_{i=0}^{M-1}a_0(i)=0$. This condition may be satisfied even for sequences $e_n$ coming from Hilbert-Kunz theory, as we show in the next example.

	\begin{example}	\label{ex:simplification of nonprimitive roots}	
		Let $\KK$ be a field of characteristic $p>0$ and let $g\geq2$. We consider the ring $R_g=\KK\llbracket s,st,\dots,st^g\rrbracket$ which is the completion of the coordinate ring of the rational normal cone of degree $g$. From \cite[Example 5.1]{Chan}, the Hilbert-Kunz function of $R_g$ with respect to its maximal ideal is
		\begin{equation*}
			\HK_{R_g}(n)=\Big(\frac{g+1}{2}\Big)p^{2n}+\frac{1}{2}(-v_n^2+v_ng-g+1),
		\end{equation*}
		where $v_n\in\{0,\dots,g-1\}$ is the congruence class of $p^n-1$ modulo $g$.
		
		When $g=5$ and $p$ is congruent to $2$ or $3$ modulo $5$, i.e., it is a generator of $(\mathbb{Z}/5\mathbb{Z})^\times$, the period of the quasipolynomial $a_0(n)=\frac{1}{2}(-v_n^2+v_ng-g+1)$ is $4$ and the sum $\sum_{i=1}^4a_0(i)$ is zero. This implies that $1$ is a root of the numerator $P(z)$ of the Hilbert-Kunz series of $\HK_{R_g}(n)$. Canceling the common factor $1-z$, the Hilbert-Kunz series takes the following form
		\begin{equation*}
			\HKS_{R_5}(z)=\frac{(3+p^2)z^3+(2+2p^2)z^2+(1+2p^2)z+1}{(1-p^2z)(1+z^2)(1+z)}.
		\end{equation*}
		When $p$ is congruent to $4$ modulo $5$, the quasi-polynomial defining the constant term of the Hilbert-Kunz series has period $2$ and no simplification occurs. The Hilbert-Kunz series is
		\begin{equation*}
			\HKS_{R_5}(z)=\frac{(-3-p^2)z^2+(1+2p^2)z+1}{(1-p^2z)(1-z)(1+z)}.
		\end{equation*}
		When $p$ is congruent to $1$ or $0$ modulo $5$, the Hilbert-Kunz function is a polynomial in $p^n$. In particular, when $p=5$, the constant term in the Hilbert-Kunz function is zero. Precisely, we have $\HK_{R_5}(n)=3\cdot 5^{2n}$
		and
		\begin{equation*}
			\HKS_{R_5}(z)=\frac{3}{1-25z}.
		\end{equation*}
		The latter fact generalizes to any $g\geq2$. If $p=g$, then the Hilbert-Kunz function is a polynomial with constant term equal to zero and the Hilbert-Kunz series is
		\begin{equation*}
			\HKS_{R_p}(z)=\frac{(p+1)/2}{1-p^2z}.
		\end{equation*}
	\end{example}
	
	\subsection{Primitive cyclotomics roots of $P(z)$}
	Now, we consider a primitive $M$-th root of unity $\zeta_M$ and we want to characterize $M$-periodic rational sequences $a_0(n)$ for which $P(\zeta_M)=0$. First, we observe that $\zeta_M$ is a root of $P(z)$ if and only if
	\begin{equation}\label{eq:radici ciclotomiche numeratore}
		P(\zeta_M)=a_0(0)-p^da_0(M-1)+\sum_{i=1}^{M-1}(a_0(i)-p^d a_0(i-1))\zeta_M^i=0.
	\end{equation}
	Since $a_0(n)\in\Q$ and $[\Q(\zeta_M):\Q]=\varphi(M)$ is the Euler totient function,  equation~\eqref{eq:radici ciclotomiche numeratore} is equivalent to a system of $\varphi(M)$ equations, linear in the variables $a_0(0),\dots,a_0(M-1)$.
 \begin{definition}
     We call $S_M$ the set of solutions of Equation~\eqref{eq:radici ciclotomiche numeratore} seen as a subspace of the ambient space $\Q^M$.
 \end{definition}
	Notice that, by construction, the $M$-uple $(a_0(0),\dots,a_0(M-1))$ lies in $S_M$ if and only if $\zeta_M$ is a root of $P(z)$. Since the linear system given by \eqref{eq:radici ciclotomiche numeratore} consists of $\varphi(M)$ equations, we have that the dimension of $S_M$ is at least $M-\varphi(M)$. We prove that it is \textit{equal} to $M-\varphi(M)$.

	\begin{theorem}\label{thm:rango massimo sistema}
		For every $M\ge 1$, we have $\dim S_M=M-\varphi(M)$.
	\end{theorem}
	\begin{proof}
		The extension $\Q(\zeta_M)/\Q$ has degree $\varphi(M)$ and $\mathcal{B}=\{1,\zeta_M,\zeta_M^2,\dots,\zeta_M^{\varphi(M)-1}\}$ is a $\Q$-basis of $\Q(\zeta_M)$.
		With respect to $\mathcal{B}$, equation \eqref{eq:radici ciclotomiche numeratore} is equivalent to the $\varphi(M)\times M$ system with matrix 
		\begin{equation*}
			\begin{pmatrix}
				1&0&0&\dots&0&0&0&-p^db_0&\dots\\
				-p^d&1&0&\dots&0&0&0&-p^db_1&\dots\\
				0&-p^d&1&\dots&0&0&0&-p^db_2&\dots\\
				\vdots&\vdots&\vdots&&\vdots&\vdots&\vdots&\vdots&\dots\\
				0&0&0&\dots&-p^d&1&0&-p^db_{\varphi(M)-3}&\dots\\
				0&0&0&\dots&0&-p^d&1&-p^db_{\varphi(M)-2}&\dots\\
				0&0&0&\dots&0&0&-p^d&1-p^db_{\varphi(M)-1}&\dots
			\end{pmatrix}
		\end{equation*}
	   in the variables $(a_0(0),a_0(1),\dots,a_0(M-1))\in\Q^M$, for suitable $b_i\in\Q$. The last written column comes mostly from the term $-p^da_0(\varphi(M)-1)\zeta_M^{\varphi(M)}$ in \eqref{eq:radici ciclotomiche numeratore}, so that the $b_i$ satisfy
		\begin{equation*}
			\zeta_M^{\varphi(M)}=b_0+b_1\zeta_M+b_2\zeta_M^2+\dots+b_{\varphi(M)-1}\zeta^{\varphi(M)-1}.
		\end{equation*}
		In other words, the $M$-th cyclotomic polynomial is
		\begin{equation*}
			\Phi_M(X)=-b_0-b_1X-b_2X^2-\dots-b_{\varphi(M)-1}X^{\varphi(M)-1}+X^{\varphi(M)}.
		\end{equation*}
		It is clear that the rank of the above matrix is at least $\varphi(M)-1$. We are going to prove that the rank is precisely $\varphi(M)$. 
		
		Let $C_0,C_1,\dots,C_{\varphi(M)-1}$ be the first $\varphi(M)$ columns of the matrix. By contradiction, write 
		\begin{equation*}
			C_{\varphi(M)-1}=\alpha_0 C_0+\alpha_1 C_1+\dots+\alpha_{\varphi(M)-2}C_{\varphi(M)-2}
		\end{equation*}
		for some $\alpha_0,\dots,\alpha_{\varphi(M)-2}\in\Q$. Looking at the first entry, we must have $\alpha_0=-b_0p^d$. Looking at the second entry, we must have $\alpha_1=-b_0p^{2d}-b_1p^d$. Proceeding in this way, see that
		\begin{equation*}
			\alpha_i=-b_0p^{(i+1)d}-b_1p^{id}-b_2p^{(i-1)d}-\dots-b_ip^{d}
		\end{equation*} 
		for every $i=0,\dots,\varphi(M)-2$. Therefore, the last entry of the linear combination gives the relation
		\begin{equation*}
			1-b_{\varphi(M)-1}p^d=-\alpha_{\varphi(M)-2}p^d,
		\end{equation*}
		that is
		\begin{equation*}
			-b_0p^{\varphi(M)d}-b_1p^{(\varphi(M)-1)d}-b_2p^{(\varphi(M)-2)d}-\dots-b_{\varphi(M)-1}p^d+1=0.
		\end{equation*}
		Dividing this relation by $p^{\varphi(M)d}$, we obtain the equation $\Phi_M(1/p^d)=0$, that is impossible since the roots of $\Phi_M(X)$ are just the primitive $M$-th roots of unity.
		This implies that the rank of the system is exactly $\varphi(M)$, hence the space $S_M$ of solutions of this system has dimension $M-\varphi(M)$.
	\end{proof}
	
	In order to give a better description of $S_M$, we introduce some notation.
	
	\begin{definition}
		For any proper divisor $\ell$ of $M$, we define 
		\begin{equation*}
			V_\ell=\{(x_0,\dots,x_{M-1})\in\Q^M \mid x_{i}=x_{\ell+i} \ \ \forall i=0,\dots, M-\ell-1 \}.	
		\end{equation*} 
	\end{definition}
	In other words, the spaces $V_\ell$ are made by all vectors of $\Q^M$ whose entries are periodic of period $\ell$. It is immediate to see that $\dim V_\ell=\ell$ and $V_{\ell_1}\cap V_{\ell_2}=V_{\gcd(\ell_1,\ell_2)}$.

	\begin{lemma}\label{lem:Vl in SM}
		For any $\ell\mid M$, $\ell\ne M$, we have that $V_\ell\subseteq S_M$.
	\end{lemma}
	\begin{proof}
		Fix $\ell\mid M$, $\ell\ne M$. A vector $(x_0,\dots,x_{M-1})\in V_\ell$ corresponds naturally to a periodic function $a(i)=x_i$ of period $\ell$. By Theorem~\ref{thm:succ-serie quasi-pol}, we know that the generating series of the function $e_n=a_dp^{nd}+a(i)$  does not have $\zeta_M$ as pole. This means that the factor corresponding to $\zeta_M$ in $Q(z)$ must simplify in $P(z)/Q(z)$. Therefore, $\zeta_M$ must be a root of $P(z)$, i.e., $(x_0,\dots,x_{M-1})\in S_M$.
	\end{proof}
	
	The lemma says that periodic functions of period $\ell$ strictly dividing $M$ are solutions of the system \eqref{eq:radici ciclotomiche numeratore}. It is natural to ask whether all solutions of \eqref{eq:radici ciclotomiche numeratore}  come from these $V_\ell$'s.
	
	\begin{question}\label{conjecture S_M}
		Is it true that $S_M=\displaystyle\sum_{\substack{\ell\mid M \\ \ell\neq M}}V_{\ell}$ holds for every $M\in\mathbb{Z}_+$?
	\end{question}

	In view of Theorem~\ref{thm:rango massimo sistema} and Lemma~\ref{lem:Vl in SM}, proving that Question~\ref{conjecture S_M} holds true is equivalent to showing that $\dim (\sum_{\ell|M, \ell\ne M}V_\ell)=M-\varphi(M)$. 
	Although the simplicity of its formulation, it seems quite hard to find a good method to give an answer in full generality. In the following, we present some partial results.
	
	\begin{proposition}
		\label{prop:case 2 prime divisors}
		Let $M=r^\alpha s^\beta$ with $r,s$ distinct odd primes, $\alpha,\beta\in\mathbb{Z}_{\ge 0}$. Then, the equality $S_M=\displaystyle\sum_{\ell|M, \ell\ne M}V_\ell$ holds.
	\end{proposition}
	\begin{proof}
		As observed above, it is enough to show that $\dim(\sum_{\ell|M, \ell\ne M}V_\ell)=M-\varphi(M)$. First, suppose $\beta=0$ (the case $\alpha=0$ is the analogous), so that $M=r^\alpha$. Then
		\begin{equation*}
			\begin{split}
				\dim\Big(\sum_{\ell|M, \ell\ne M}V_\ell\Big)&=\dim V_{r^{\alpha-1}}=r^{\alpha-1}=r^{\alpha}-(r^\alpha-r^{\alpha-1})=M-\varphi(M).
			\end{split}
		\end{equation*}

		If both $\alpha$ and $\beta$ are greater than zero, we have that
		\begin{equation*}
			\begin{split}
				\dim\Big(\sum_{\ell|M, \ell\ne M}V_\ell\Big)&=\dim(V_{r^\alpha s^{\beta-1}}+V_{r^{\alpha-1}s^\beta})=\dim(V_{r^\alpha s^{\beta-1}})+\dim(V_{r^{\alpha-1}s^\beta})-\dim(V_{r^{\alpha-1}s^{\beta-1}})=\\
				&=r^\alpha s^{\beta-1}+r^{\alpha-1}s^\beta-r^{\alpha-1}s^{\beta-1}=r^\alpha s^{\beta}-(r^\alpha-r^{\alpha-1})(s^\beta-s^{\beta-1})=\\
				&=M-\varphi(M).
			\end{split}
		\end{equation*}
	\end{proof}
		
        Again, let $e_n=a_d p^{nd}+a_0(n)$ be a rational sequence and call $M$ the period of $a_0(n)$. Summing up the main result of this section, we have found that
        \begin{itemize}
            \item[(a)] If $a_0(n)$ is a linear combination of sequences of period strictly smaller than $M$, then $\zeta_M$ is a root of $P(z)$ (Lemma \ref{lem:Vl in SM});
            \item[(b)] If $M$ is a power of a prime number, $\zeta_M$ is a root of $P(z)$ if and only if $a_0(n)$ has exact period strictly smaller than $M$ (Proposition \ref{prop:case 2 prime divisors});
            \item[(c)] If $M$ has exactly two distinct prime factors, $\zeta_M$ is a root of $P(z)$ if and only if $a_0$ is a linear combination of sequences of period strictly smaller than $M$ (Proposition \ref{prop:case 2 prime divisors}). 
        \end{itemize}
 
	Generalizing the argument of the proof of Proposition \ref{prop:case 2 prime divisors} to the case where $M$ has three or more prime factors seems complicated due to the lack of a Grassman-like formula for the sum of more than two subspaces. However, we carried on some explicit examples and all of them seem to verify the content of Question \ref{conjecture S_M}.
	
	When $M$ is divisible by more than one prime, one can then produce many examples of sequences $e_n=a_d p^{nd}+a_0(n)$ where $M$ is the exact period of $a_0(n)$ and $\zeta_M$ is a root of $P(z)$. A more subtle question is whether there exist Hilbert-Kunz functions with this property. The answer is positive, as shown by the following example.

	\begin{example}\label{ex:rootofunitysimplifies}
		Let $p=2$ and let $\KK$ be a field of characteristic $p$. We consider the rings $R_g=\KK\llbracket s,st,\dots,st^g\rrbracket$ from Example~\ref{ex:simplification of nonprimitive roots} for $g=3$ and $g=7$ .
		We can realize both rings as modules over the power series ring $S=\KK\llbracket x_1,\dots,x_8\rrbracket$ in $8$ variables as follows. For $R_7$, we take the scalar multiplication given by the map $\varphi:S\rightarrow R_7$ such that $\varphi(x_i)=st^{i-1}$. For $R_3$, we take the scalar multiplication given by the map  $\psi:S\rightarrow R_3$ such that $\psi(x_i)=st^{i-1}$ if $i\leq4$ and $\psi(x_i)=0$ if $i\geq5$.
		Let $\mathfrak{m}=(x_1,\dots,x_8)$ be the maximal ideal of $S$. Notice that the maximal ideal of $R_7$ coincides with $\mathfrak{m}R_7$, and the same is true for $R_3$. This implies that the Hilbert-Kunz functions of $R_3$ and $R_7$ as $S$-modules (with respect to the maximal ideal $\mathfrak{m}$) coincide with the Hilbert-Kunz functions of $R_3$ and $R_7$ as local rings with respect to their maximal ideals.
		By \cite[Example 5.1]{Chan}, these functions are given by
		\begin{equation*}
			\begin{split}
				\HK_{R_3}(n)=\begin{cases}
					2\cdot 2^{2n}-1 &\text{for $n$ even}\\
					2\cdot 2^{2n} &\text{for $n$ odd}
				\end{cases}
			\end{split},
		\end{equation*}
		with constant term of period $2$, and by
		\begin{equation*}
			\begin{split}
				\HK_{R_7}(n)=\begin{cases}
					4\cdot 2^{2n}-3 &\text{for $n\equiv 0\pmod 3$}\\
					4\cdot 2^{2n} &\text{for $n\equiv 1\pmod 3$}\\
					4\cdot 2^{2n}+3 &\text{for $n\equiv 2\pmod 3$}
				\end{cases}
			\end{split},
		\end{equation*}
		with constant term of period $3$.
		
		The sum of these rational sequences gives a sequence whose constant term has exact period $6$. The Hilbert-Kunz function of the $S$-module $R_7\oplus R_3$ is precisely this sequence, that is  
  \[
  \HK_S(R_3\oplus R_7;n):=\dim_\KK\left(R_3\oplus R_7/\mathfrak{m}^{[p^n]}(R_3\oplus R_7)\right)=6\cdot 2^{2n}+a_0(n),
  \] 
  where $a_0(n)$ takes periodically the values $-4,0,2,-3,-1,3$. By construction, $a_0(n)$ is the sum of two periodic functions of period respectively $2$ and $3$, therefore, by Lemma~\ref{lem:Vl in SM}, the numerator $P(z)$ of the Hilbert-Kunz series $\HKS_S(R_7\oplus R_3;z)$ has $\zeta_6$ as one of its roots. After simplifying the primitive factor $z^2-z+1$, we obtain
		\begin{equation*}
			\HKS_S(R_7\oplus R_3;z)=\frac{18z^4+11z^3-18z^2-18z-2}{(1-4z)(z-1)(z+1)(z^2+z+1)}.
		\end{equation*}
	\end{example}

	

\bibliographystyle{siam}
\bibliography{references}

\end{document}